\documentclass[11pt]{article}
\usepackage{epsfig}
\usepackage{multirow}
\newcommand{\be}{\begin{equation}}
\newcommand{\ee}{\end{equation}}

\renewcommand{\theequation}{\arabic{section}.\arabic{equation}}

\begin{document}

\newcommand{\ra}{\rightarrow}
\newcommand{\wh}{\widehat}
\newcommand{\nit}{\noindent}
\newcommand{\no}{\nonumber}
\newcommand{\ba}{\begin{eqnarray}}
\newcommand{\ea}{\end{eqnarray}}
\newcommand{\pa}{\mbox{$\partial$}}
\newcommand{\Gam}{\mbox{$\Gamma$}}
\newcommand{\dta}{\mbox{$\delta$}}
\newcommand{\fee}{\mbox{$\varphi$}}
\newcommand{\al}{\mbox{$\alpha$}}
\newcommand{\lam}{\mbox{$\lambda$}}
\newcommand{\eps}{\mbox{$\epsilon$}}
\newcommand{\gam}{\mbox{$\gamma$}}
\newcommand{\omg}{\mbox{$\omega$}}
\newtheorem{theo}{Theorem}[section]
\newtheorem{defin}{Definition}[section]
\newtheorem{prop}{Proposition}[section]
\newtheorem{lem}{Lemma}[section]
\newtheorem{cor}{Corollary}[section]
\newtheorem{rmk}{Remark}[section]
\renewcommand{\theequation}{\arabic{section}.\arabic{equation}}

\title{A Dynamic Algorithm for Blind Separation of Convolutive Sound Mixtures}

\author{Jie ${\rm Liu}^{*}$, \hspace{.1 in} Jack Xin\thanks{Department of Mathematics,
UC Irvine, Irvine, CA 92697, USA.}\hspace{.05 in},
 \hspace{.02 in} and
\hspace{.02 in} Yingyong Qi \thanks{
Qualcomm Inc,
5775 Morehouse Drive,
San Diego, CA 92121, USA.}
}
\date{}
\thispagestyle{empty}

\maketitle

\thispagestyle{empty}

\begin{abstract}
We study an efficient dynamic blind source separation algorithm of convolutive
sound mixtures based on updating statistical information in the frequency domain, and
minimizing the support of time domain demixing filters by a weighted least square method.
The permutation and scaling indeterminacies of separation, and
concatenations of signals in
adjacent time frames are resolved with optimization of
$l^1 \times l^\infty$ norm on cross-correlation coefficients
at multiple time lags. The algorithm is a direct method without iterations, and
is adaptive to the environment. Computations on recorded and
synthetic mixtures of speech and music signals show excellent performance.


\end{abstract}

\vspace{.2 in}

\hspace{.2 in} {\it Keywords}: Convolutive Mixtures,
Indeterminacies,

\hspace{1.0 in}  Dynamic Statistics Update, Optimization,

\hspace{1.0 in} Blind Separation.

\newpage
\setcounter{page}{1}
\section{Introduction}
\setcounter{equation}{0}
Blind source separation (BSS) methods aim to extract
the original source signals from their mixtures based on the statistical
independence of the source signals without knowledge of the mixing environment.
The approach has been very successful for instantaneous mixtures. However, realistic
sound signals are often mixed through a media channel, so
the received sound mixtures are linear convolutions of the unknown sources and the channel transmission functions.
In simple terms, the observed signals are unknown
weighted sums of the signals and its delays. Separating convolutive mixtures
is a challenging problem especially in realistic settings.

In this paper, we study a dynamic BSS method using both frequency and time domain
information of sound signals in addition to the
independence assumption on source signals.
First, the convolutive mixture in the time domain is decomposed
into instantaneous mixtures in the frequency domain by the fast
Fourier transform (FFT). At each frequency, the joint approximate
diagonalization of eigen-matrices (JADE) method is applied.
The JADE method collects second and fourth order statistics from
segments of sound signals to form a set of matrices for joint
orthogonal diagonalization, which leads to an estimate of de-mixing
matrix and independent sources. However, there remain extra degrees
of freedom: permutation and scaling of estimated sources at each frequency.
A proper choice of these parameters is critical for the separation
quality. Moreover, the large number of samples of the statistical approach can
cause delays in processing.  These issues are to be addressed by utilizing
dynamical information of signals in an optimization framework.
We propose to dynamically update statistics with newly received signal frames,
then use such statistics to determine permutation in the frequency domain by
optimizing an $l^1 \times l^\infty$ norm of
channel to channel cross-correlation coefficients with multiple time lags.
Though cross channel correlation functions and related similarity measure
were proposed previously to fix permutation \cite{Murata01}, they allow cancellations
and may not measure similarity as accurately and reliably as the norm (metric) we introduced here.
The freedom in scaling is fixed by minimizing the support of the estimated
de-mixing matrix elements in the time domain. An efficient weighted least square method
is formulated to achieve this purpose directly in contrast to iterative method
in \cite{PS00}. The resulting dynamic
BSS algorithm is both direct and adapted to the acoustic environment. Encouraging results on
satisfactory separation of recorded sound mixtures are reported.

The paper is organized as follows. In section
2, a review is presented on frequency domain approach, cumulants and joint
diagonalization problems, and indeterminacies. Then the proposed
dynamic method is presented, where objective functions of optimization,
statistics update and efficient
computations are addressed. Numerical results are shown and analyzed to demonstrate
the capability of the algorithm to separate speech and music mixtures
in both real room and synthetic environments. Conclusions are in section 3.

\section{Convolutive Mixture and BSS}
Let a real discrete time signal be
$s(k) = [s_1(k),s_2(k),\cdots,s_{n}(k)]$, $k$ a discrete time index,
such that the components $s_i(k)$ ($i=1,2,\cdots,n$), are
zero-mean and mutually independent random processes. For simplicity,
the processing will divide $s$ into partially overlapping frames of length $T$ each.
The independent components are transmitted and mixed to give the
observations $x_i(k)$:
\be
x_i(k) = \sum_{j=1}^{n}\,  \sum_{p=0}^{P-1}\, a_{ij}(p)\, s_{j}(k -p), \; i=1,2,\cdots, n; \label{bs1}
\ee
where $a_{ij}(p)$ denote mixing filter coefficients, the $p$-th element of the $P$-point
impulse response from source $i$ to receiver $j$. The mixture in (\ref{bs1}) is
convolutive, and an additive Gaussian noise may be added. The sound signals
we are interested in are
speech and music, both are non-Gaussian \cite{BS95}.
We shall consider the case of equal number of
receivers and sources, especially $n=2$.

An efficient way to decompose the nonlocal equation (\ref{bs1}) into local ones is by
a T-point discrete Fourier transform (DFT) \cite{Br},
$X_j(\omega,t)= \sum_{\tau=0}^{T-1}\, x_{j}(t+\tau)\,e^{-2\pi \,J\omega \tau}$,
where $J=\sqrt{-1}$, $\omega$ is a frequency index, $\omega = 0, 1/T, \cdots, (T-1)/T$,
$t$ the frame index.
Suppose $T >P$, and extend $a_{ij}(p)$ to all $p \in [0, T-1]$ by zero padding.
Let $H_{ij}(\omega)$ denote the matrix function obtained by T-point DFT of $a_{ij}(p)$ in $p$,
$S_j(\omega,t)$ the T-point DFT of $s_j(k)$ in the $t$-th frame.
If $P \ll T$, then to a good approximation \cite{PS00}:
\be
X(\omega,t)\approx H(\omega)\, S(\omega, t), \label{bs2}
\ee
where $X=[X_1,\cdots,X_n]^{Tr}$, $S=[s_1,\cdots,s_n]^{Tr}$, $Tr$ is short for transpose.
The components of $S$ remains independent of
each other, the problem is converted to a blind separation of instantaneous
mixture in (\ref{bs2}). Note that $P$ is on the order of 40 to 50 typically,
while $T$ is 256 or 512, so the assumption $P\ll T$ is reasonable.

\subsection{Instantaneous Mixture and JADE}
Let us briefly review an efficient and accurate method,
so called joint approximate diagonalization eigen-matrices (JADE)
\cite{Car93} for BSS of instanteneous mixture.  There are many other approaches in the
literature \cite{Choi05},
e.g. info-max method \cite{BS95} which is iterative and based on maximizing some information
theoretical function. JADE is essentially a direct method for reducing covariance.
We shall think of $S$ as a random function of $t$, and suppress $\omega$ dependence.
First assume that by proper scaling $E[|S_j(t)|^2]=1$, $j=1,\cdots,n$.
It follows from independence of
sources that ($'$ conjugate transpose):
\be
E[ S(t)S(t)^{'}] = I_{n}, \;\; R_X \equiv E[ X(t)X(t)^{'}] = H H^{'}, \label{bs3}
\ee
the latter identity is a factorization of the Hermitian covariance matrix of the mixture.
However, there is non-uniqueness in the ordering and phases of columns of $H$.
Suppose that (1) the mixing matrix $H$ is full rank; (2) the $S_j(t)$'s are independent at any $t$;
(3) the process $S(t)$ is stationary. Let $W$ be a matrix
such that $I_n = W R_X W^{'} = WH H' W'$, $W$ is called a whitening matrix. Then
$WH$ is an orthogonal matrix, denoted by $U$. Multiplying $W$ from the left onto (\ref{bs2}),
one finds that:
\be
Z(t)\equiv W X(t) = U S(t). \label{bs4}
\ee
The 4th order statistics are needed to determine $U$. The 4th-order cumulant of
four mean zero random variables is:
\be
{\rm Cum}[a,b,c,d]=E(abcd) - E(ab)E(cd) -E(ac)E(bd)-E(ad)E(bc), \label{bs5}
\ee
which is zero if $a,b,c,d$ split into two mutually independent groups. For source vector $S$,
Cum$[S_i,S_j,S_k,S_l]= {\rm kurt}_i\, \delta_{ijkl}$, ${\rm kurt}_{i}$ = Cum$[S_i,S_i,S_i,S_i]$
is the kurtosis. If ${\rm kurt}_{i} \not = 0$, the i-th source is called kurtic. Kurtosis is zero
for a mean zero Gaussian random variable.
The last assumption of JADE is that (4) there is at most one non-kurtic source.

Define cumulant matrix set $Q_{Z}(M)$ from $Z$ in (\ref{bs4}) as the linear span of the
Hermitian matrices $Q=(q_{ij})$ satisfying ($*$ complex conjugate):
\be
q_{ij} = \sum_{k,l=1}^{n}\, {\rm Cum}(Z_i,Z_{j}^{*},Z_{k},Z_{l}^{*})\, m_{lk},\;\; 1\leq i,j\leq n,
\label{bs6}
\ee
where matrix $M=(m_{ij})=e_{l}e_{k}^{'}$, $e_l$ being the unit vector with zero components except
the $l$-th component equal to one. Equations (\ref{bs4}) and (\ref{bs6})
imply that ($u_p$ is the p-th column of $U$):
\be
Q= \sum_{p=1}^{n}\, ({\rm kurt}_{p}\, u_{p}^{'} M u_{p})\, u_{p}u_{p}^{'}, \;\; \forall \; M,
\label{bs7}
\ee
or $Q = U D U^{'}$, $D= {\rm diag}({\rm kurt}_{1}u_{1}^{'}Mu_{1},\cdots,{\rm kurt}_{n}u_{n}^{'}Mu_{n})$.
Hence, $U$ is the joint diagonalizer of the matrix set $Q_{Z}(M)$.  Once $U$ is so determined,
the mixing matrix $H=W^{-1} U$. It can be shown \cite{Car93} using identity (\ref{bs7}) that the
joint diagonalizer of $Q_Z(M)$ is
equal to $U$ up to permutation and phase, or up to a matrix multiplier $P$ where $P$
has exactly one unit modulus entry in each row and column. Such a joint diagonalizer is
called essentially equal to $U$.

The algorithm of finding the joint diagonalizer is a generalization of Jacobi method or
Givens rotation method \cite{GV}. As the cumulant matrices are estimated in practice,
exact joint diagonalizer may not exist, instead, an approximate joint diagonalizer, an
orthogonal matrix $V$, is sought
to maximize the quantity:
$C(V,B) = \sum_{r=1}^{n^2}\, |{\rm diag}(V'\, B_r\, V)|^{2}$,
where $B=\{B_1,B_2,\cdots, B_{n^2}\}$ is a set of basis (or eigen)
matrices of $Q_Z(M)$, $|{\rm diag}(A)|^2$ is the sum of squares of diagonals of a matrix $A$.
Maximizing $C(V,B)$ is same as minimizing off diagonal entries, which can be achieved in
a finite number of steps of Givens rotations. The costs of joint diagonalization is
roughly $n^2$ times that of diagonalizing a single Hermitian matrix.

Though stationarity is assumed for the theoretical analysis above, JADE turns out to be
quite robust even when stationarity is not exactly satisfied for signals such as speech
or music.

\subsection{Dynamic Method of Separating Convolutive Mixture}
For each frequency $\omega$, equation (\ref{bs2}) is a BSS problem of instantaneous mixtures.
The speech or music signals in reality are stationary over short time scales and
nonstationary over longer time scales, which depend on the production details.
For speech signals, human voice is stationary for a few 10 ms, and becomes non-stationary
for a time scale above 100 ms due to envelope modulations \cite{dengo,Murata01}.
The short time stationarity permits FFT to generate meaningful spectra in equation
(\ref{bs2}) within each frame. For
a sampling frequency of 16,000 Hertz, each frame of 512 points lasts 32 ms.
The mixing matrix $H$ may depend on $t$ over longer time scales,
denoted by $H=H(\omega,t)$,
unless the acoustic environment does not change as in most synthetic mixing.
A demixing method with potential real time application should be able to
capture the dynamic variation of mixing.

Our approach consists of four steps.
Step I is to find an initialization for $H(\omega,t)$. After receiving the initial $n_T$
frames of mixtures, compute their FFT and obtain $X(\omega,t)$, $t=1,2,\cdots, n_T$, to
collect $n_T$ samples at each discrete frequency. For each $\omega$, perform JADE,
and estimate the mixing matrix denoted by $H_0(\omega)$. To ensure a good statistical
estimate, $n_T$ is on the order of $80$ to $100$, and may be properly reduced later.

Step I gives separated components of signals over all frequencies.
However, such JADE output has inderterminacies in amplitude, order and phase.
This benign problem for instantaneous mixtures becomes a major issue when one needs to
assemble the separated individual components. For example, the
permutation mismatches across frequencies can degrade the quality of separation seriously.

Step II is to use nonstationarity of signals to sort out a consistent order of
separated signals in the frequency domain. Such a method for batch processing was
proposed in \cite{Murata01}. A separation method requiring the entire length of the
signal is called batch processing.
The sorting algorithm of \cite{Murata01} proceeds as follows. (1) Estimate the envelope variation
by a moving average over a number of frames (beyond stationarity time scale)
for each separated frequency component.
The envelope is denoted by ${\rm Env}(\omega,t,i)$, where $i$ is the index
of separated components. (2) Compute a
similarity measure equal to the sum of correlations of the envelopes
of the separated components at each frequency. The similarity measure is
 ${\rm sim}(\omega) = \sum_{i\not = j}\, \rho({\rm Env}(\omega,t,i),{\rm Env}(\omega,t,j))$,
where $\rho(\cdot,\cdot)$ is the normalized correlation coefficients (see (\ref{bs11})) involving
time average over the {\it entire signal
length} to approximate the ensemble average so the $t$ dependence drops out. (3) Let $\omega_1$ be
the one with lowest similarity value where separation is the best. The $\omega_1$ serves as
a reference point for sorting. (4)
At other frequencies $\omega_k$ ($k =2,3,\cdots$),
find a permutation $\sigma$ to maximize $
 \sum_{i=1}^{n}\, \rho({\rm Env}(\omega_k,t,\sigma(i)),
\sum_{j=1}^{k-1} {\rm Env}_s(\omega_j,t,i))$,
among all permutations of $1,2,\cdots,n$. Here ${\rm Env_s}$ denotes the sorted
envelopes in previous frequencies. (5) Permute the order of separated components
at the $k$-th frequency bin according to $\sigma $ in step (4), and
define ${\rm Env}_s(\omega_k,t,i)$. Repeat (4) and (5) until $k=T$.

We shall modify the above sorting method in three aspects. The first is to
use segments of signal instead
of the entire signal to compute statistics (correlations) to
minimize delay in processing.
The second is to use correlation coefficients of
separated signals at {\it un-equal times} or {\it multiple time lags}
in step (2) to better characterize the degree of
separation. Moreover, we notice that the similarity measure of \cite{Murata01} as seen
above is a sum of correlation coefficients of potentially both signs, and so
can be nearly zero due to cancellations even though each term in the sum is not small in absolute value.
We introduce an $l^1\times l^\infty$ norm below to characterize more accurately
channel similarity by taking sum of absolute values of correlation coefficients and
maximum of time lags. The third is to simplify the maximization problem on $\sigma$
to avoid comparing correlations with summed envelopes at all previous frequencies. 
We also do not use envelopes of signals inside correlation functions. The reason is that 
the smoothing nature of envelope operation reduces the amount of oscillations in the signals
and may yield correlation values less accurate for capturing the degree of independence.
Specifically, let $\hat{s}_{i}(\omega,t)=a_{i}(\omega,t)\, e^{j \phi_{i}(\omega,t)}$
be the $i$-th separated signal at frequency $\omega$, where
$a_{i}(\omega,t) = \left|\hat{s}_{i}(\omega,t)\right|$,
$\phi_i$ the phase functions, $t$ the frame index.
The correlation function of two time dependent signals over
$M$ frames is:
\be
{\rm cov}(a(\omega,t),b(\omega',t))= M^{-1} \sum_{t=1}^{M}\, a(\omega,t)b^{*}(\omega',t) -
M^{-2}\sum_{t=1}^{M} a(\omega,t) \, \sum_{t=1}^{M} b^{*}(\omega',t), \label{bs10}
\ee
and the (normalized) correlation coefficient is:
\be
\rho(a(\omega,t),b(\omega',t))={{\rm cov}(a(\omega,t),b(\omega',t)) \over
\sqrt{{\rm cov}(a(\omega,t),a(\omega,t))\, {\rm cov}(b(\omega',t),b(\omega',t))}}.
\label{bs11}
\ee
From speech production viewpoint, frequency components of a speech signal do not
change drastically in time, instead are similarly affected by
the motion of the speaker's vocal chords. The correlation coefficient is
a natural tool for estimating coherence of frequency components of a speech signal. A similar argument may be applied to
music signals as they are produced from cavities of instruments.

Now with $M =n_T$ in (\ref{bs10}), define
\begin{equation} \label{C-omg}
  C(\omega)= \sum_{i\not = j} \,\max_{k\in \{-K_0,...,K_0\} } \,
  |\rho(|\hat{s}_{i}(\omega,t)|,|\hat{s}_{j}(\omega, t -k)|)|, \quad \mbox{for } \omg\in[\omg_L,\omg_U]
\end{equation}
with some positive integer $K_0$.
Find $\omega_1$ between $\omg_L$ and $\omg_U$ to minimize
$C(\omega)$. With $\omega_1$ as reference, at any other
$\omega$, find the permutation $\sigma $ to maximize:
\begin{equation} \label{f-permutation}
  \sigma ={\rm argmax}\, \sum_{i=1}^{n}\,\max_{k\in \{-K_0,...,K_0\}}\,
  |\rho(|\hat{s}_{i}(\omega_1,t)|, |\hat{s}_{\sigma(i)}(\omega,t-k)|)|.
\end{equation}
Notice that the objective functions in (\ref{C-omg})-(\ref{f-permutation})
are exactly the $l^{1} \times l^{\infty}$ norms over the indices $i$($j$) and $k$.
Multiple time lag index $k$ is to accomodate the translational invariance of
sound quality to the ear. Maximizing over $k$ helps to capture the correlation of the
channels, and sum of $i$ ($j$) reflects the total coherence of a vector signal.

Step III fixes the scaling and phase indeterminacies in $\hat{s}(\omega,t)$.
Each {\em row} of the de-mixing matrix $H^{-1}_{0}(\omega)$ may be multiplied by
a complex number $\lam_i(\omega)$ ($i=1,2\cdots,n$) before inverse FFT (ifft) to reconstruct demixing
matrix $h^{(0)}(\tau)$ in the time domain. The idea is to minimize the support of each row of
the inverse FFT by a weighted least square method. In other words,
we shall select $\lam_i$'s so that the
entries of ${\rm ifft}(H^{-1}_{0})(\tau)\equiv h^{(0)}(\tau)$
are real and nearly zero if $\tau \geq Q$ for some $Q < T$, $Q$ as small as possible,
$T$ being the length of FFT. Smaller $Q$ improves the local approximation, or accuracy of
equation (\ref{bs2}). To be more specific, using $H^{-1}_{0,i}(\omega)$
to denote the $i$-th {\em row vector} of $H^{-1}_{0}(\omega)$, we can
explicitly write the equation to shorten the support of inverse FFT:
\be \label{scaling-1}
  \mbox{ifft}(\lam_i(\omega)  H^{-1}_{0,i}(\omega))(\tau) = 0
\ee
in terms of the real and imaginary parts of $\lam_i(\omg)$ for $\omg=0,1/T,...,(T-1)/T$.
Those real and imaginary parts are the variables and the equations are {\em linear}.
Now, we let $\tau$ run from $q$ to $T-1$.
If we want small support, $q$ should be small, then there are more equations
than unkowns. So we multiply a weight to each equation and minimize
in the least square sense. Equation (\ref{scaling-1}) for larger $\tau$ is
multiplied by a larger weight
in the hope that the value of the left hand side of (\ref{scaling-1})
will be closer to zero during the least square process.
If we choose the weighting function to be the exponential function $\beta^\tau$
for some $\beta>1$, then the above process can be mathematically written as
\be  \label{bs12}
  [\lam_{i}(0),...,\lam_i((T-1)/T)] = \mbox{argmin}
   \sum_{\tau=q}^{T-1} |\beta^\tau \mbox{ifft}(\lam_i(\omega)  H^{-1}_{0,i}(\omega))(\tau)|^2
\ee
where $H^{-1}_{0,i}(\omega)$ is the $i$-th {\em row vector} of $H^{-1}_{0}(\omega)$.

A few comments are in order.
First, since the mixing matrix $H_0(\omg)$ is the FFT of
a real matrix, we impose that $H_0(\omega)=H_0(1-\omega)^*$. So,
supposing $T$ is even, we only need to apply JADE to obtain
$H_0(\omg)$ for $\omg=0,1/T,...,1/2$; $H_0(0)$ and $H_0(1/2)$ will automatically
be real. When fixing the freedom of
scaling in each $\omega$, we choose $\lam(0)$ and $\lam(1/2)$ real,
and $\lam(\omg)=\lam(1-\omg)^*$ for other $\omg$. Second, to fix the overall scaling and
render the solution nontrivial, we set $\lam(0) = 1$.
Third, the weighted least square problem (\ref{bs12}) can be solved by a
direct method or matrix inversion (chapter 6 in \cite{GV}).

Note that when $n=2$, among the $2(T-q)$ equations from
(\ref{scaling-1}) with $\tau=q,...,T-1$,
there are $T-1$ variables including $\lam_i(1/2)$, the real and imaginary parts of
$\lam_i(\omg)$ for $\omg=1/T,...,1/2-1/T$. So, we can make roughly
half of $h_i^{(0)}(\tau) \approx 0$, the best one can achieve in general.
Separated signals, denoted by $\tilde{s}^{(0)}(t)$,
are then produced, for $t \in [0,n_T]$, $t$ the frame index.

The last step IV is to update $h^{(0)}(\tau)$ when $\delta n_T \ll n_T$ many new frames of
mixtures arrive. The steps I to III are repeated using frames from $\delta n_T + 1$ to
$\delta n_T + n_T$, to generate a new time domain demixing matrix $h^{(1)}(\tau)$, $\tau\in[0,T-1]$, and
separated signal $\tilde{s}^{(1)}(\tau)$, $\tau \in [T(n_T- \Delta n_T)+1 ,T(n_T + \delta n_T)]$
with $T$ the size of one frame. We use $\tau$ here instead of $t$ because in the most part of the paper,
$t$ is the frame index. Now, $\tilde{s}^{(1)}(\tau)$ and $s^{(0)}(\tau)$ share a
common interval of size $T \Delta n_T$.
On this common interval, $\tilde{s}^{(1)}(\tau)$ and $s^{(0)}(\tau)$ will
be the same if we are doing a perfect job and if the ordering of $\tilde{s}^{(1)}$ is consistent
with that of $s^{(0)}$. In order to determine the ordering of $\tilde{s}^{(1)}(\tau)$,
we compute $\rho\left(\tilde{s}^{(0)}_i (\tau), s_j^{(1)}(\tau-k)\right)$
on this common interval with different $k$ and $i,j=1,...,n$.
Then we determine the permutation $\sigma$ of the components of $\tilde{s}^{(1)}(t)$
by minimization:
\begin{equation} \label{t-permutation}
   \sigma = {\rm argmax}\, \sum_{i=1}^{n} \max_{k\in\{-K_1,...,K_1 \}} \,
   \left|\rho\left( s_{i}^{(0)}(\tau), \tilde{s}^{(1)}_{\sigma(i)}(\tau-k)\right)\right|
\end{equation}
with some constant $K_1$. After doing the necessary permutation of $\tilde{s}^{(1)}$, the separated signals
are then extended to the extra frames $\delta n_T + n_T$
by concatenating the newly separated $\delta n_T$ many frames of $\tilde{s}^{(1)}$ with
those of $\tilde{s}^{(0)}$. The continuity of concatenation is maintained
by requiring that $\max_{\tau}\, |h_{ii}^{(k)}(\tau)|$'s ($i=1,2,\cdots,n$) are invariant
in $k$, where $k=1,2,\cdots$, labels the updated filter matrix in time.
 The procedure repeats with the next arrival of mixture data, and
is a direct method incorporating dynamic information.
\medskip

\noindent Because sorting order depends only on the
relative values of channel correlations, we observed in practice that
the $\max_{k\in \{-K_.,...,K_.\}}$ in
equations (\ref{C-omg}), (\ref{f-permutation}),
(\ref{t-permutation}) may be replaced by
$\sum_{k=-K_.}^{K_.}$, with a different choice of $K_.$ value. The
$\max_{k\in \{-K_.,...,K_.\}}$ is a more accurate characterization however.

\medskip

\subsection{Adaptive Estimation and Cost Reductions}
Cumulants and moments are symmetric functions
in their arguments \cite{NP93}.
For example when $n=2$, there are 16 joint fourth order cumulants from (\ref{bs5}), however,
only six of them need to be computed, the others follow from symmetry. Specifically,
among the 16 cumulants:
\[
  \begin{array}{rr}
    Q(1)={\rm Cum}(y_1,y^*_1,y^*_1,y_1), &  Q(2)={\rm Cum}(y_1,y^*_1,y^*_1,y_2)\\
    Q(3)={\rm Cum}(y_1,y^*_1,y^*_2,y_1), &  Q(4)={\rm Cum}(y_1,y^*_1,y^*_2,y_2)\\
    Q(5)={\rm Cum}(y_1,y^*_2,y^*_1,y_1), &  Q(6)={\rm Cum}(y_1,y^*_2,y^*_1,y_2)\\
    Q(7)={\rm Cum}(y_1,y^*_2,y^*_2,y_1), &  Q(8)={\rm Cum}(y_1,y^*_2,y^*_2,y_2)\\
    Q(9)={\rm Cum}(y_2,y^*_1,y^*_1,y_1), &  Q(10)={\rm Cum}(y_2,y^*_1,y^*_1,y_2)\\
    Q(11)={\rm Cum}(y_2,y^*_1,y^*_2,y_1), &  Q(12)={\rm Cum}(y_2,y^*_1,y^*_2,y_2)\\
    Q(13)={\rm Cum}(y_2,y^*_2,y^*_1,y_1), &  Q(14)={\rm Cum}(y_2,y^*_2,y^*_1,y_2)\\
    Q(15)={\rm Cum}(y_2,y^*_2,y^*_2,y_1), &  Q(16)={\rm Cum}(y_2,y^*_2,y^*_2,y_2)
  \end{array}
\]
we have the relations: $Q(2)=Q(3)^*=Q(5)^*=Q(9)$, $Q(4)=Q(6)=Q(11)=Q(13)$,
 $Q(7)=Q(10)^*$, $Q(8)=Q(15)=Q(12)^*=Q(14)^*$, where $*$ is complex conjugate.
For N samples, we only need to compute the
following six $1\times N$ vectors
\[
 \begin{array}{rr}
  Y_1=(y^1_1 y^1_1,...,y^N_1 y^N_1), &
  Y_2=(y^1_1 y^1_2,...,y^N_1 y^N_2), \\
  Y_3=(y^1_2 y^1_2,...,y^N_2 y^N_2), &
  Y_4=(y^1_1 y^{1*}_1,...,y^N_1 y^{N*}_1), \\
  Y_5=(y^1_1 y^{1*}_2,...,y^N_1 y^{N*}_2), &
  Y_6=(y^1_2 y^{1*}_2,...,y^N_2 y^{N*}_2), \\
 \end{array}
\]
then all the 4th order and 2nd order statistical quantities can be reconstructed. For example,
\begin{equation} \label{Q1}
  Q(1)= \frac{1}{N}Y_4 \cdot Y_{4}^{Tr}  - \frac{1}{N^2}\left( 2 \; {\rm sum}(Y_4) \; {\rm sum}(Y_4) + {\rm sum}(Y_1) \;
  \left({\rm sum}(Y_1)^*\right)\right)
\end{equation}
where ${\rm sum}(Y_i)$ is the summation of the $N$ components of
$Y_i$.

As formula (\ref{bs5}) suggests, cumulants are updated through moments
when $\delta n_T$ early samples are replaced by
the same number of new samples. As $\delta n_T$ is much less than the total number of
terms $n_T$ in the empirical estimator of expectation, the adjustment costs $2 \delta n_T$ flops
for each second moments and $6 \delta n_T$ flops for each joint fourth order moment.
The contributions of the
early samples are subtracted from the second and fourth moments, then
the contributions of the new samples are added. The cumulant update approach is
similar to cumulant tracking method of moving targets (\cite{LMen99} and references therein).

Due to dynamical cumulants update, the prewhitening step at each frequency is
performed after cumulants are computed from $X(\omega)$.
This is different from JADE \cite{Car93} where the prewhitening occurs before computing the commulants.
This way, it is more convenient to make use of the previous cumulant information and
updated $X(\omega)$. Afterward, we use the
multilinearity of the cummulants to transform them back to the commulants of
the prewhitened $X(\omega)$, before joint diagonalization.

It is desirable to decrease $n_T$ to lower the number of samples for cumulants estimation.
However, this tends to increase the variance in the estimated cumulants, and render estimation
less stable in time. Numerical experiments
indicated that with $n_T$ as low as 40, the separation using overlapping frames is still
reliable with reasonable quality.

It is known \cite{dengo} that the identity of a speaker is carried by
pitch (perception of the fundamental frequency in speech production)
which varies in the low frequency range of a few hundred Hertz. We found that
instead of searching among all frequencies for the reference frequency $\omega_1$ in step II,
it is often sufficient to search in the low frequency range. The smaller searching range
alleviates the workload in sorting and permutation correcting. This is similar to a feature oriented
method, see \cite{Osher_90,qxin,chan98} among others.

\subsection{Experimental Results}
The proposed algorithm with adaptivity and cost reduction considerations
was implemented in Matlab. The original code of JADE by J.-F. Cardoso
is obtained from a open source
(http://web.media.mit.edu/$\sim$paris/ ) maintained by P. Smaragdis.
Separation results with both dynamic and batch processing of three different
types of mixtures are reported here:
\begin{itemize}
 \item[(1)] real room recorded data;
 \item[(2)] synthetic mixture of speech and music;
 \item[(3)] synthetic mixture of speech and speech noise.
\end{itemize}
They will be called case (1), (2) and (3) in the following discussion.

The values of the parameters used in the three cases are listed in Table~\ref{T.parameters}.
In the table, "$n_T$ (dyn.)" is the initial value of $n_T$ in dynamic processing
and "$n_T$ (bat.)" is the $n_T$ in batch processing. Other than $n_T$, dynamic and batch
process share the same parameters. The frame size is $T$,
"overlap" is the overlapping percentage between two successive frames,
$\delta n_T$ and $\Delta n_T$ are as in step IV, $K_0$ and $K_1$ are from
(\ref{f-permutation}) and (\ref{t-permutation}), $\beta$ is in (\ref{bs12}), and
$q$ is the lower limit of $\tau$ in (\ref{scaling-1}).

Note that the values of $\omg_L$ and $\omg_U$ from (\ref{C-omg})
are not listed in the table. In our computation, we use the
following two choices
\begin{itemize}
 \item[(A)] $\omg_L=0$, $\omg_U=1/2$.
 \item[(B)] $\omg_L=\omg_U=4/T$, namely fixing reference frequency $\omg_1=4/T$.
\end{itemize}
For the three cases reported in this paper, {\em both choices work
and generate very similar results.} As a consequence, we will only
plot the results of the first choice. The first choice is more
general while the second is motivated by the pitch range of speech
signal and is computationally more favorable. However, we do not know precisely
the robustness of the latter.

\begin{table}[hbtp]
    \begin{center}
        \begin{tabular}{|l|l|l|l|l|l|l|l|l|l|l|l|l|l|l|}
            \hline
             case   &  T  & overlap & $n_T$ & $\delta n_T$ & $\Delta n_T$
                    &  $K_0$  & $K_1$ &  $\beta$  & $q$ &  $n_T$ \\
             & & & (dyn.) &&&&&& & (bat.) \\ \hline

            (1)      &  512 & 0\% & 100    & 20  & 30  &  4 & 10 & 1.04 & 2 & 200 \\\hline
            (2)      &  256 & 50\% & 100   & 20  & 40  & 15 & 20 & 1.04 & 2 & 160 \\\hline
            (3)      &  256 & 50\% & 100   & 20  & 40  & 10 & 20 & 1.04 & 2 & 160 \\\hline
        \end{tabular}
    \end{center}
    \caption{Parameters used in both dynamic and batch processing.}
\label{T.parameters}
\end{table}

For a quantitative measure of separation in all three cases, we compute
the maximal correlation coefficient over multiple time lags:
\begin{equation}
   \bar{\rho}(a,b)=\max_{k \in \{ -K_2,...,K_2 \} } \left| \rho(a(\tau),b(\tau+k)) \right|
\end{equation}
with $\rho$ defined in (\ref{bs11}). The $\bar{\rho}$ is computed for the mixtures,
the sources and the separated signals for both batch and dynamic
processing. An exception is the lack of sources in case (1). We
choose $K_2=20$ in all the computations. The results are listed in
Table~\ref{T.rho} which shows that the $\bar{\rho}$ values of the mixtures 
are much larger than those of the dynamically separated signals, which are on the same order as 
the $\bar{\rho}$ values of the batch separated signals. In the synthetic cases (2) and (3), 
the $\bar{\rho}$ values of the batch separated signals are on the same order of the 
$\bar{\rho}$ values of the source signals or $10^{-2}$. In cases (2) and (3), we use the 
ratio $\bar{\rho}(x,s_1) / \bar{\rho}(x,s_2)$ to measure the relative closeness of a signal 
$x$ to source signals $s_1$ and $s_2$. Table~\ref{T.rho2} lists these ratios for $x$ being the 
separated signals by dynamic and batch methods with $A$ and $B$ denoting the two ways of setting the 
reference frequency $\omega_1$. The outcomes are similar no matter $x=\tilde{s}_1$ or $x=\tilde{s}_2$ 
(first or second separated signal) in either dynamic or batch cases and either way of selecting the reference frequency $\omega_1$.

\begin{table}[hbtp]
    \begin{center}
        \begin{tabular}{|c|c|c|c|c|c|c|c|c|c|}
            \hline
             $\bar{\rho}(\cdot,\cdot)$ of 3 cases   &  mixture & dyn. separation  & bat. separation & sources \\ \hline
            (1)-A   & 0.8230 & 0.0269  & 0.0160  & N/A \\
            (1)-B   & 0.8230  & 0.0225  & 0.0159  & N/A\\\hline

            (2)-A   & 0.6240  & 0.0503 & 0.0673  & 0.0201  \\
            (2)-B   & 0.6240  & 0.0182 & 0.0600  & 0.0201 \\ \hline
            (3)-A   & 0.4613  & 0.0351  & 0.0378  & 0.0243 \\
            (3)-B   & 0.4613 & 0.0267  & 0.0677  & 0.0243 \\ \hline

        \end{tabular}
    \end{center}
    \caption{Values of the correlation coefficient $\bar{\rho}(x,y)$, $(x,y)$ being either the two  
mixtures or the two sources or the two separated signals by dynamic and batch methods.
The A and B in the first column denote the two different ways of
selecting the reference frequency $\omg_1$. }
\label{T.rho}
\end{table}

\begin{table}[hbtp]
\begin{center} {
\begin{tabular}{|r|r|r|}
\hline
      
$\bar{\rho}(x,s_1)/\bar{\rho}(x,s_2)$
      &  case(2) &
        case(3) \\ \hline
x= dyn. $\tilde{s}_1(A)$   &  4.5899 &  4.5096 \\
x= dyn. $\tilde{s}_2(A)$   &  0.1086 &  0.2852  \\ \hline
x= dyn. $\tilde{s}_1(B)$   &  5.3083 &   5.8411 \\
x= dyn. $\tilde{s}_2(B)$   & 0.0494  & 0.2799    \\ \hline 
x= bat. $\tilde{s}_1(A)$   &  15.0912 &  1.4632 \\
x= bat. $\tilde{s}_2(A)$   & 0.0760 & 0.1665 \\ \hline
x= bat. $\tilde{s}_1(B)$   & 6.2227  &  25.8122\\
x= bat. $\tilde{s}_2(B)$   & 0.0636  & 0.1719 \\ \hline
\end{tabular} }
\end{center}
\caption{Ratios of $\bar{\rho}(x,s_1)$ and $\bar{\rho}(x,s_2)$, $x$ being 
a separated signal on the first column by dynamic or batch method, 
$s_1$ and $s_2$ are source signals. The ratio measures the relative 
closeness of $x$ to $s_1$ and $s_2$. If the ratio is larger (smaller) than one, 
$x$ is closer to $s_1$ ($s_2$).
The A and B in the first column denote the two different ways of
selecting $\omg_1$. } \label{T.rho2}
\end{table}

In case (1), the recorded data \cite{Murata01} consists of 2 mixtures of a piece of
music (source 1) and a digit (one to ten) counting
sentence (source 2) recorded in a normal office size room.
The sampling frequency is 16 kHz, and 100 k data points are shown in Fig. 1.
The signals last a little over 6 seconds.  The result of dynamic BSS algorithm is shown
in Fig 2. As a comparison, we show in Fig. 3 result of batch processing of steps I to III of
the algorithm with $n_T=200$. The batch processing gives a clear separation upon
listening to the separated signals. The dynamic processing is comparable.
The filter coefficients in the time domain $h_{ij}(\tau)$
at the last update of dynamic processing are shown in Fig. 4.
Due to weighted least square
optimization in step III, they are localized and
oscillatory with support length $Q$ close to half of the FFT size $T$.

For cases (2) and (3), we show the {\em envelopes of the absolute values} of the
mixtures or the separated signals. The signal envelope was computed using the standard
procedure of amplitude demodulation, i.e., lowpass
filtering the rectified signal. The filter was an FIR
filter with 400 taps and the cutoff frequency was 100
Hz. Signal envelopes help to visualize and compare
source and processed signals.
We have normalized all the envelopes
so that the maximum height is 1. The values of $a_{ij}$ in (\ref{bs1}), which are
used to synthetically generate the mixtures, are shown in Fig.~\ref{aij} (see \cite[(8)]{Tok96}).
Fig.~\ref{f-m.1} and Fig.~\ref{f-m.2} show the mixtures and
separated signals of case (2).
Fig.~\ref{f-n.1} and Fig.~\ref{f-n.2} show the mixtures and
separated signals of case (3). In view of these plots, Table 2 and Table 3,
separation is quite satisfactory, which is also confirmed by hearing the separated
signals.

The processing time in MATLAB on a laptop can be a factor of 5 to 8 above the real time
signal duration, however, the time is expected to be closer to real time with
the computation is executed by Fortran or C directly or with additional
cost reduction techniques.
A breakdown of time consumption in the algorithm shows that
40\% of the processing is spent on computing cumulants, 30 \% on
sorting in frequency and time domains, 15\% on fixing scaling functions, 3\% on
joint diagonalization, the rest on other operations such as computing lower order statistics,
FFT, IFFT etc.

\section{Conclusions}
A dynamic blind source separation algorithm is proposed to track the
time dependence of signal statistics and to be adaptive to the potentially
time varying environment.
Besides an efficient updating of cumulants, the method made precise the
procedure of sorting permutation indeterminacy in the
frequency domain by optimizing a metric (the $l^1\times l^\infty$ norm) on
multiple time lagged channel correlation coefficients. A direct and efficient
weighted least square approach is introduced to compactify the support of
demixing filter to improve the
accuracy of frequency domain localization of convolutive mixtures.
Experimental results show robust and satisfactory separation of real recorded data and
synthetic mixtures. An interesting line of future work will be concerned with various
strategies to reduce computational costs.

\section{Acknowledgements}
The work was partially supported by NSF grants ITR-0219004, DMS-0549215, NIH grant 2R44DC006734;
the CORCLR (Academic Senate Council on Research, Computing and Library Resources) faculty
research grant MI-2006-07-6, and a Pilot award of the Center for Hearing Research at UC Irvine.


\newpage
\begin{center}
Figure Captions
\end{center}

\noindent Fig 1: Case (1), two recorded signals in a real room where a
speaker was counting ten digits with music playing in the
background.

\medskip

\noindent Fig 2: Case (1) with choice A, separated digit counting sentence
(bottom) and background music (top) by the proposed dynamic
method. Choice B gives similar results.

\medskip

\noindent Fig 3: Case (1) with choice A, separated digit counting sentence
(bottom) and background music (top) by batch processing using the
proposed steps I to III. Choice B gives similar results.

\medskip

\noindent Fig 4: Case (1) with choice A, the localized and oscillatory
filter coefficients in the time domain at the last frame of
dynamic processing. Choice B gives similar results.

\medskip

\noindent Fig 5: The weights $a_{ij}$ used in generating synthetic mixtures
of cases (2) and (3), as proposed in \cite{Tok96}.

\medskip

\noindent Fig 6: Case (2), the synthetic mixtures are generated by a female
voice and a piece of instrumental music.

\medskip

\noindent Fig 7: Case (2) with choice A, the envelopes of the separated
signals from mixtures whose envelopes are in Fig. (\ref{f-m.1}).
The small amplitude portion of the music is well recovered. Choice
B gives similar results.

\medskip

\noindent Fig 8: Case (3), the synthetic mixtures of a female voice and a
speech noise with signal to noise ratio equal to $-3.8206$ dB. The
$x_1$ plot shows a speech in a strong noise, the valley structures
in the speech signal are filled by noise.

\medskip

\noindent Fig 9: Case (3) with choice A, the envelopes of the separated
signals, noise (top) and speech (bottom). The envelopes of the two
mixtures are in Fig.~\ref{f-n.1}. The strongly noisy $x_1$ in
Fig.~\ref{f-n.1} has been cleaned, the valleys in the envelope
re-appeared. Choice B gives an even better result.

\newpage

\begin{figure}
\centerline{\includegraphics[width=400pt]{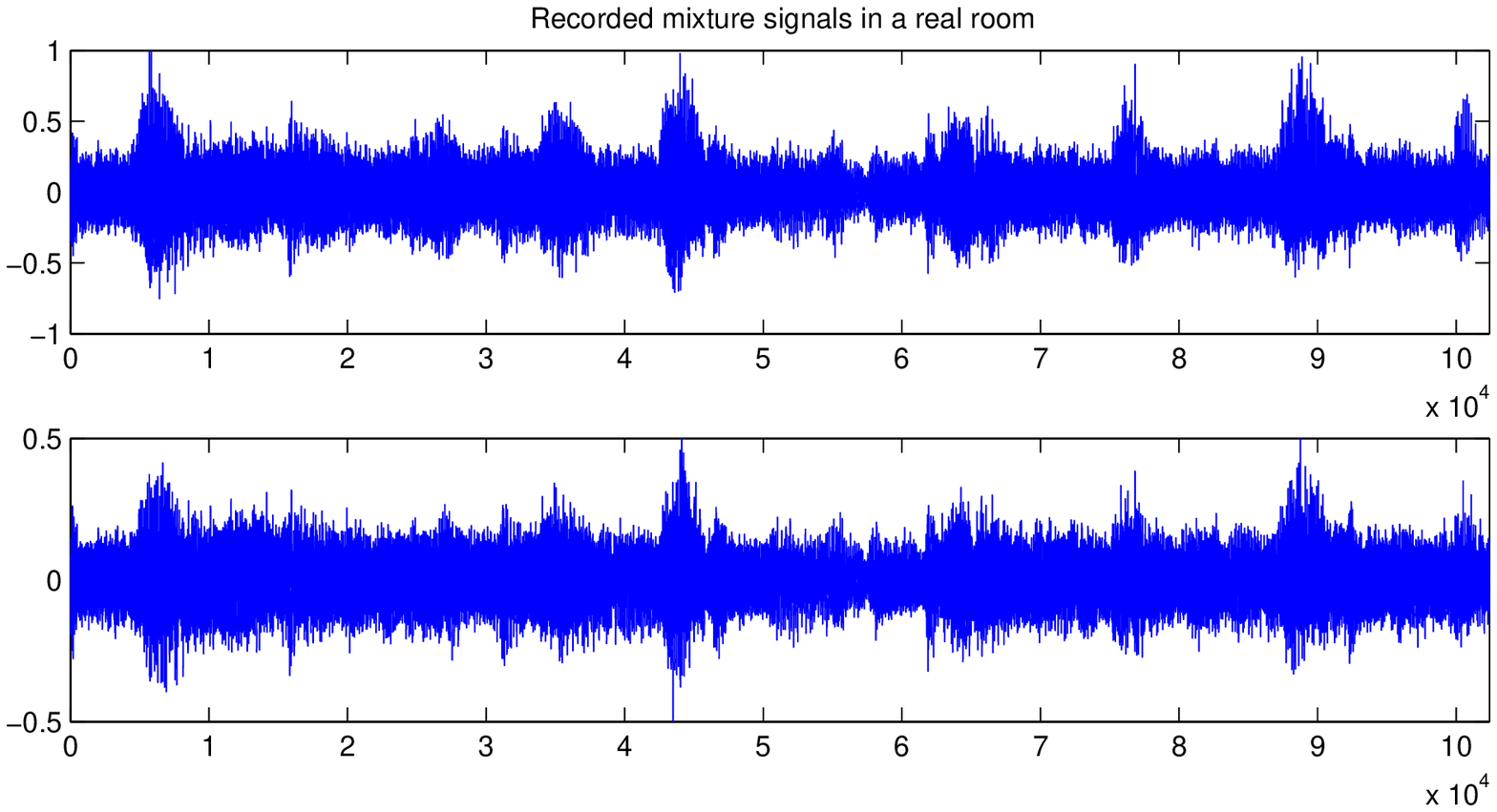}}
\caption{ }
\end{figure}

\newpage

\begin{figure}
\centerline{\includegraphics[width=400pt]{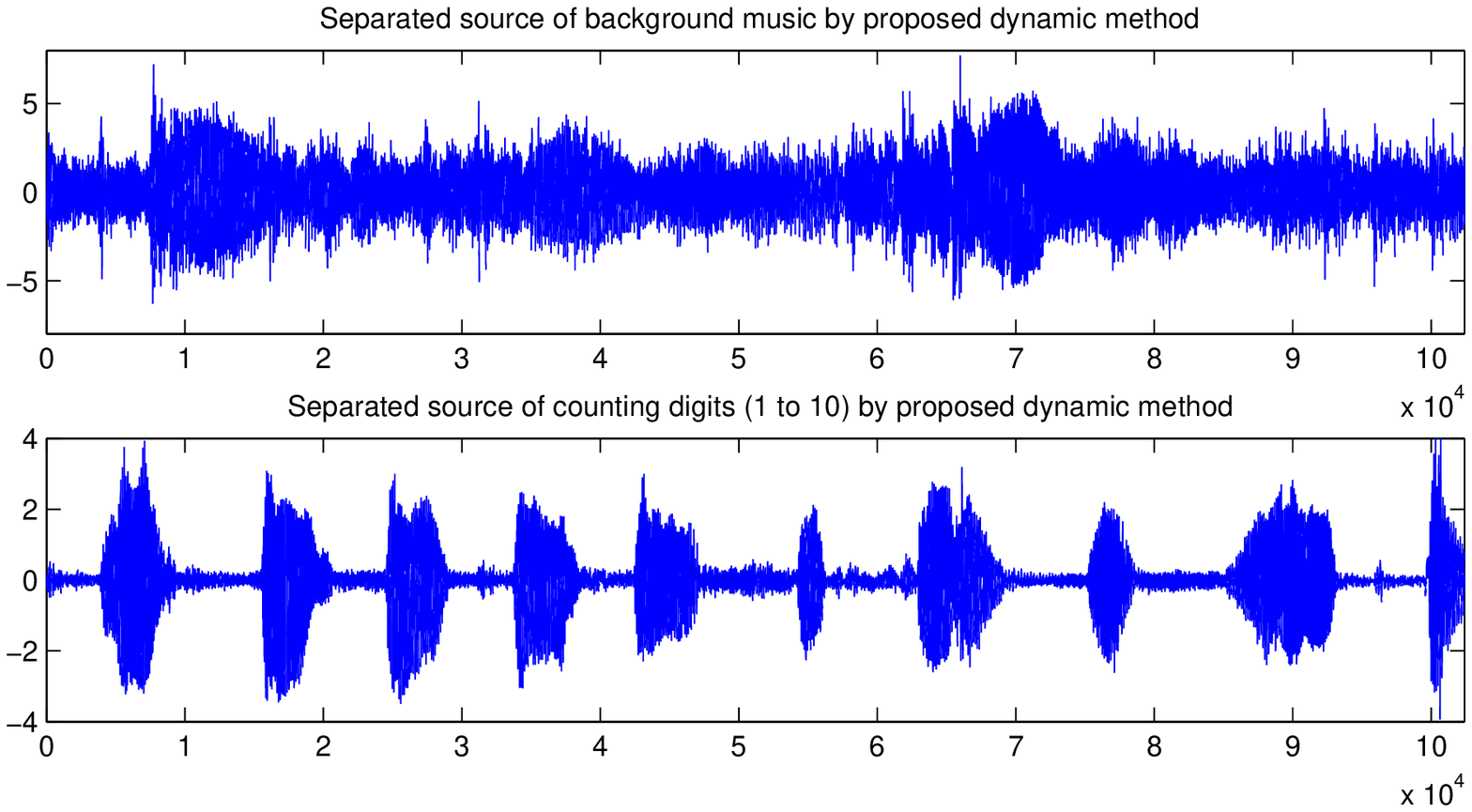}}
\caption{}
\end{figure}

\newpage

\begin{figure}
\centerline{\includegraphics[width=400pt]{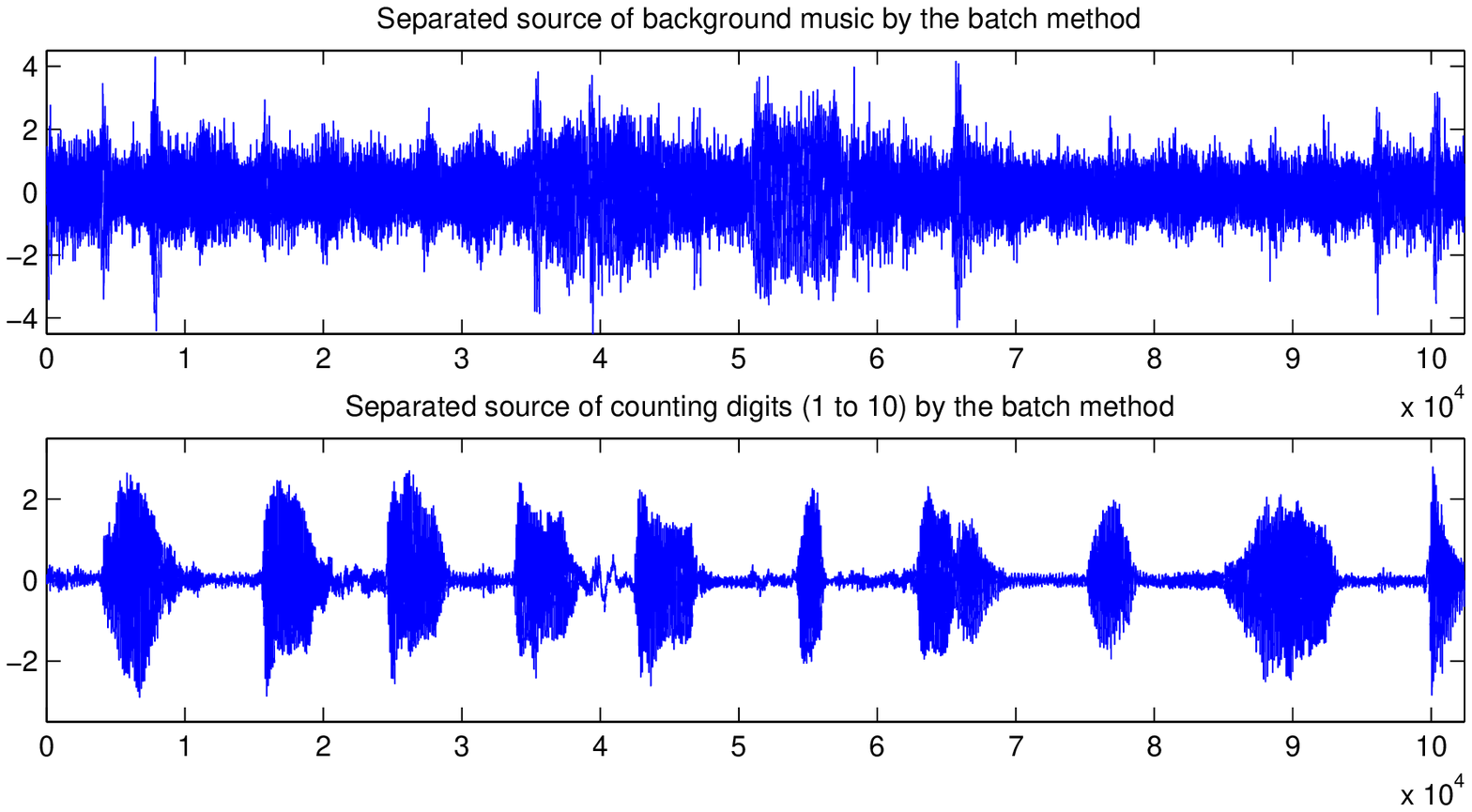}}
\caption{}
\end{figure}

\newpage

\begin{figure}
\centerline{\includegraphics[width=400pt]{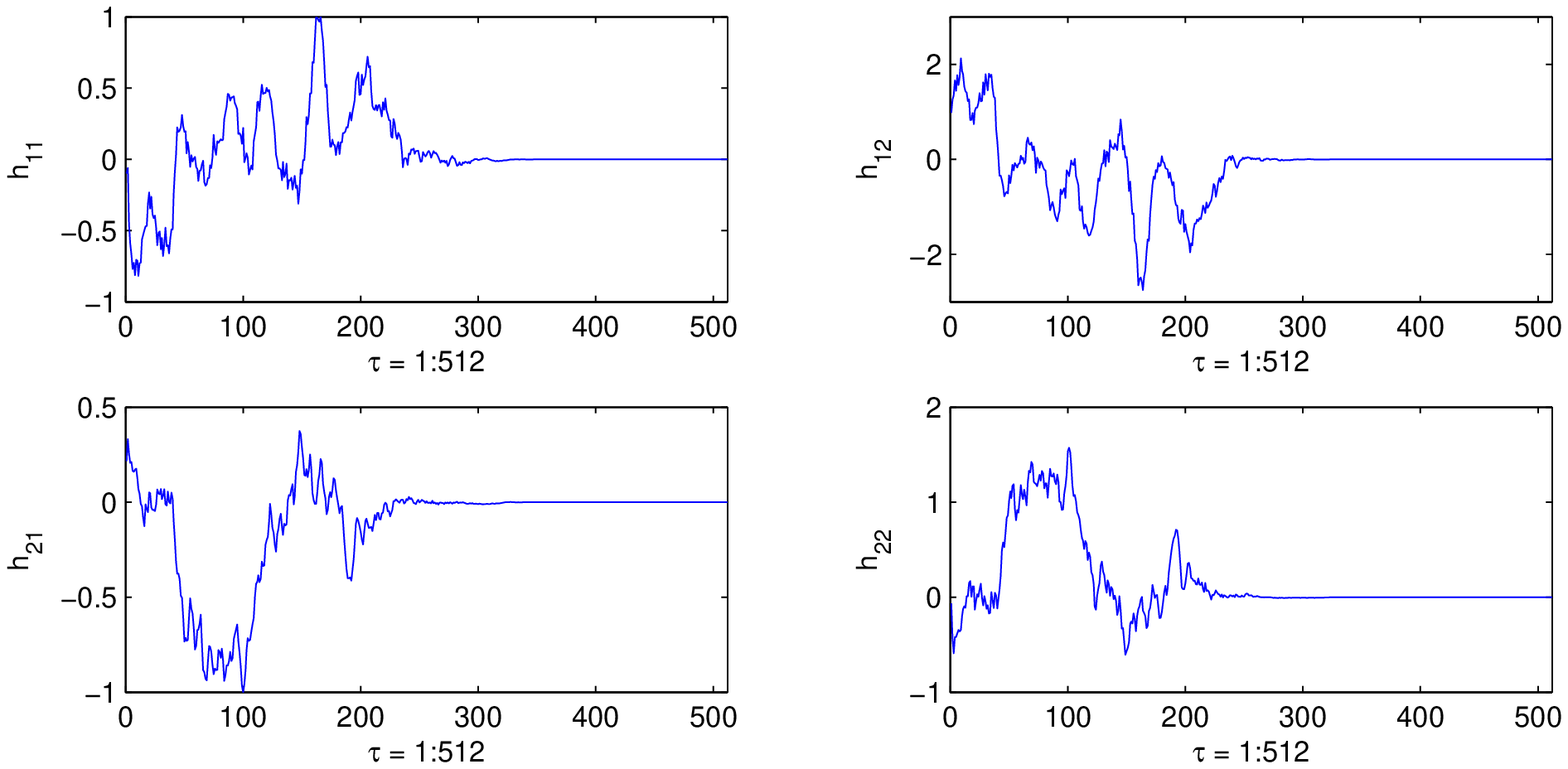}}
\caption{}
\end{figure}

\newpage

\begin{figure}
\centerline{\includegraphics[width=400pt]{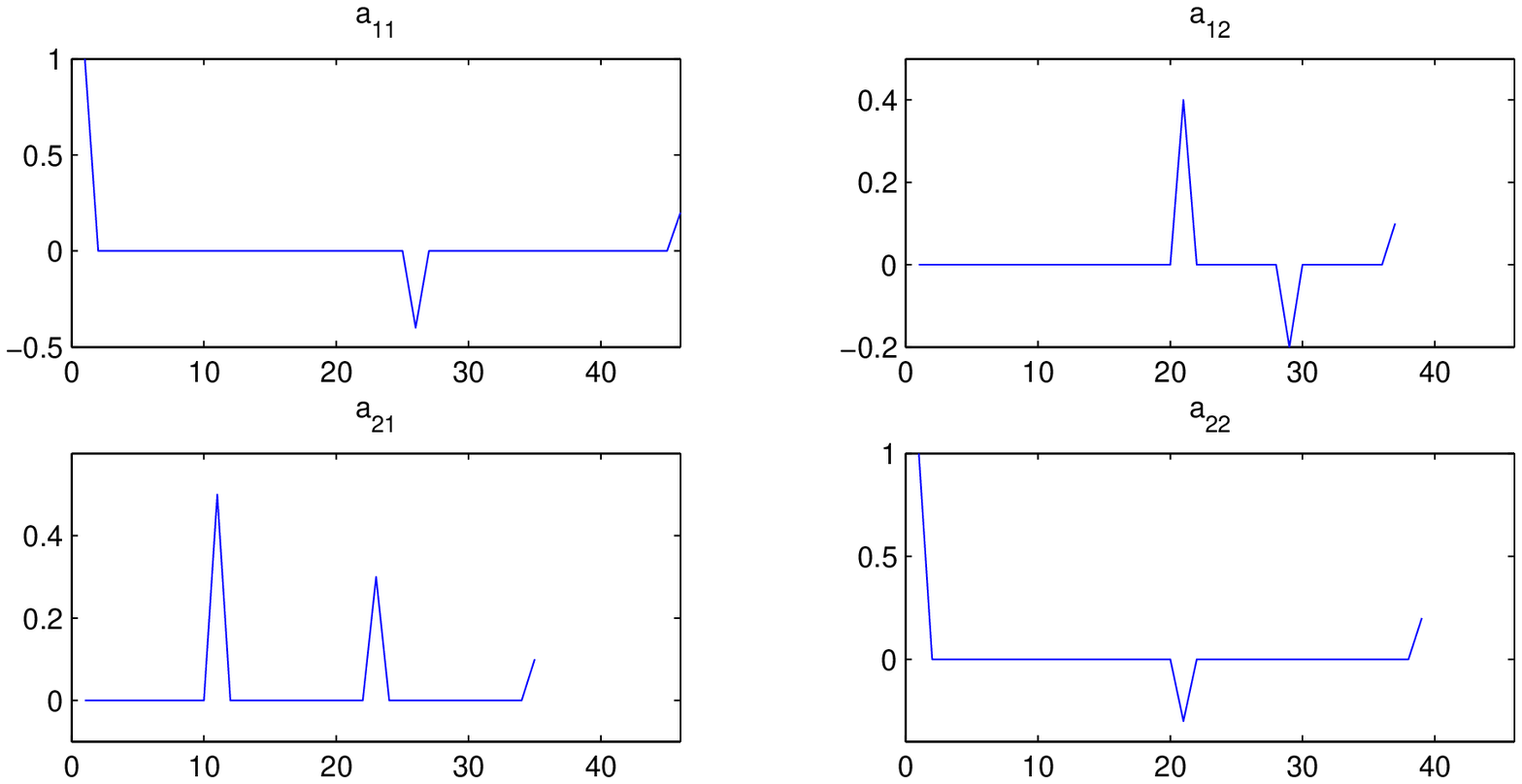}}
\caption{} \label{aij}
\end{figure}

\newpage

\begin{figure}
\centerline{\includegraphics[width=400pt]{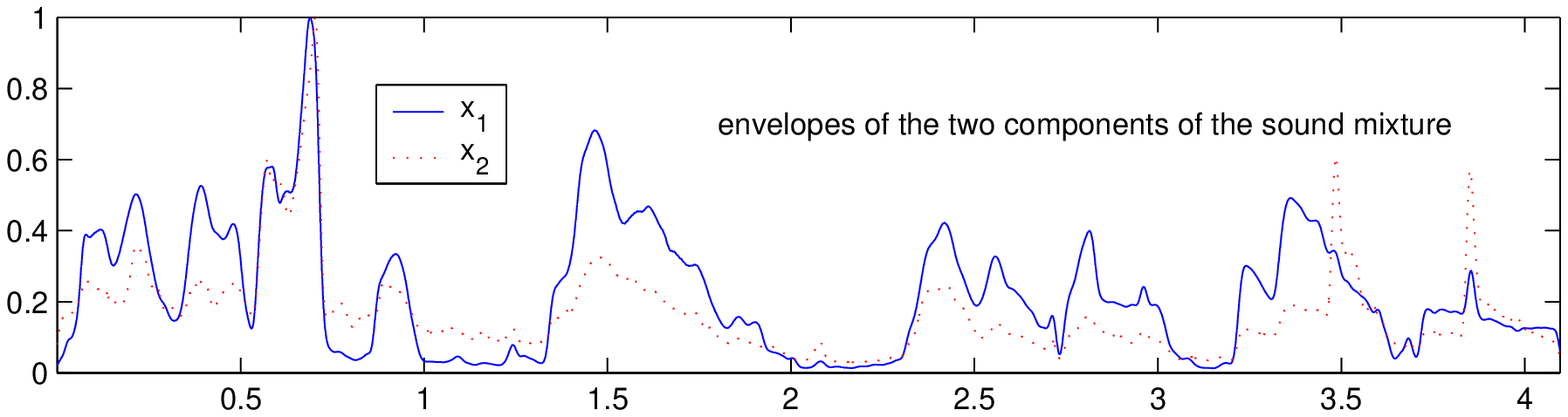}}
\caption{} \label{f-m.1}
\end{figure}

\newpage

\begin{figure}
\centerline{\includegraphics[width=400pt]{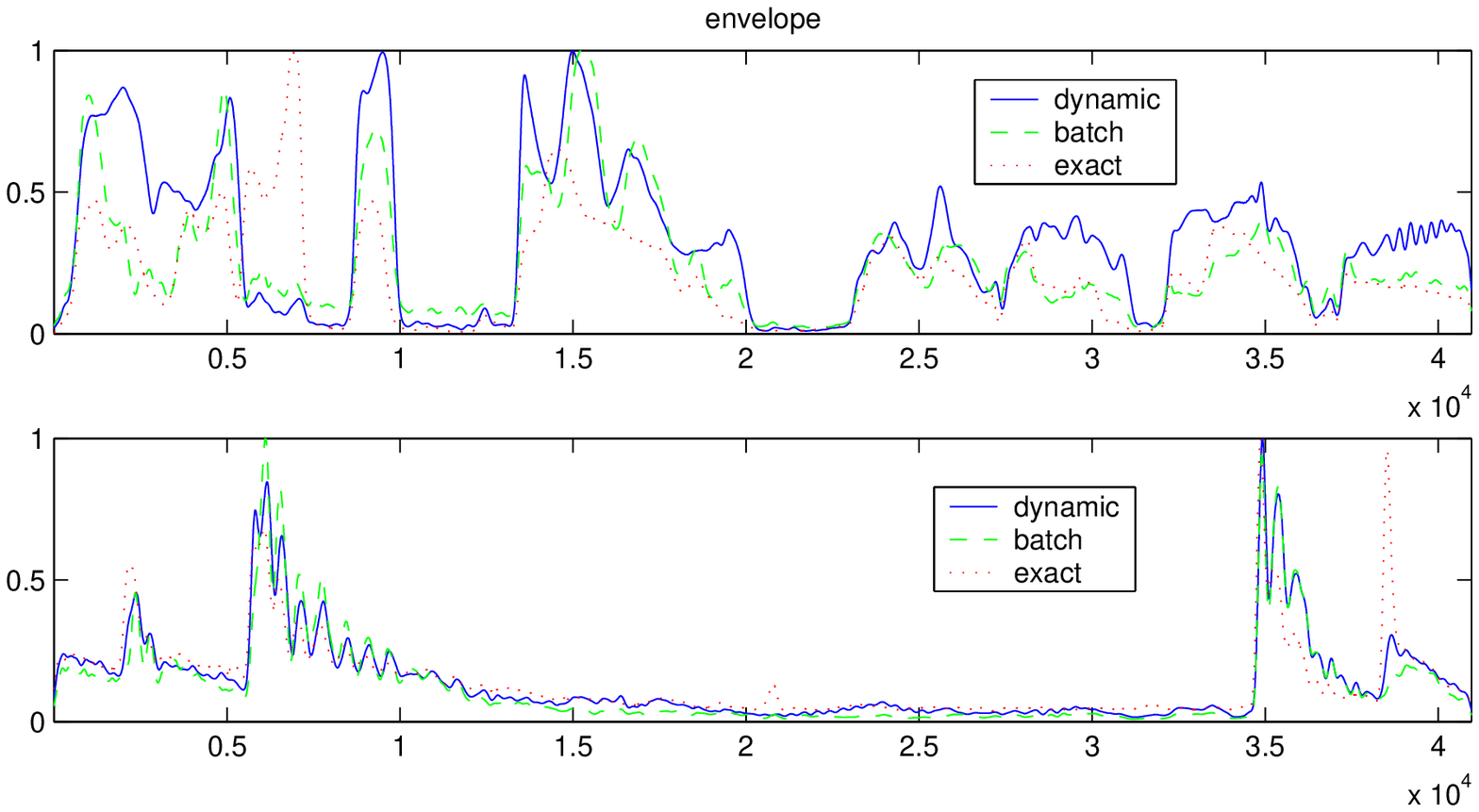}}
\caption{} \label{f-m.2}
\end{figure}

\newpage

\begin{figure}
\centerline{\includegraphics[width=400pt]{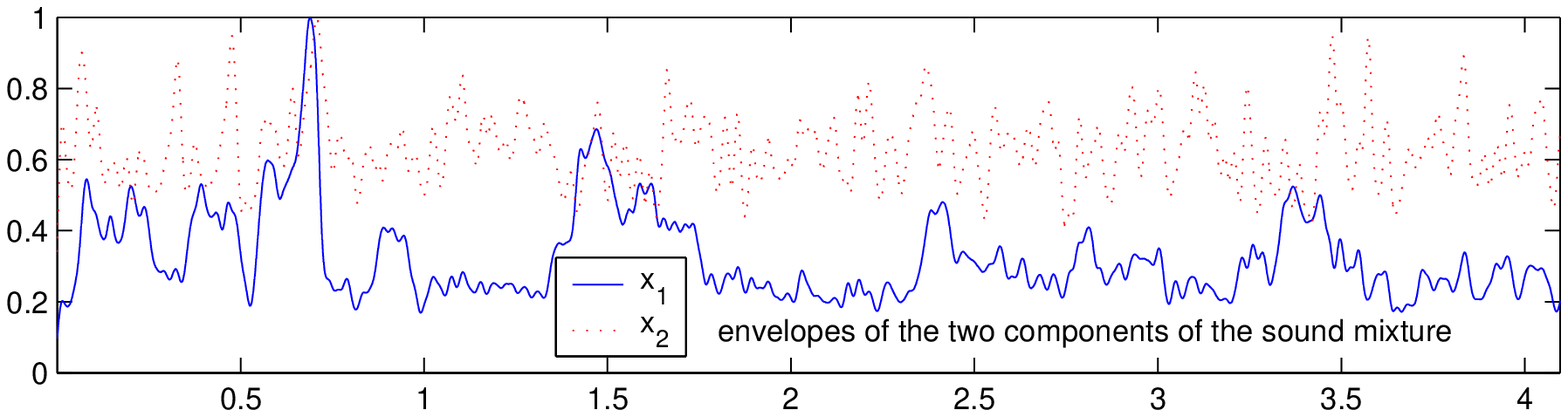}}
\caption{} \label{f-n.1}
\end{figure}

\newpage

\begin{figure}
\centerline{\includegraphics[width=400pt]{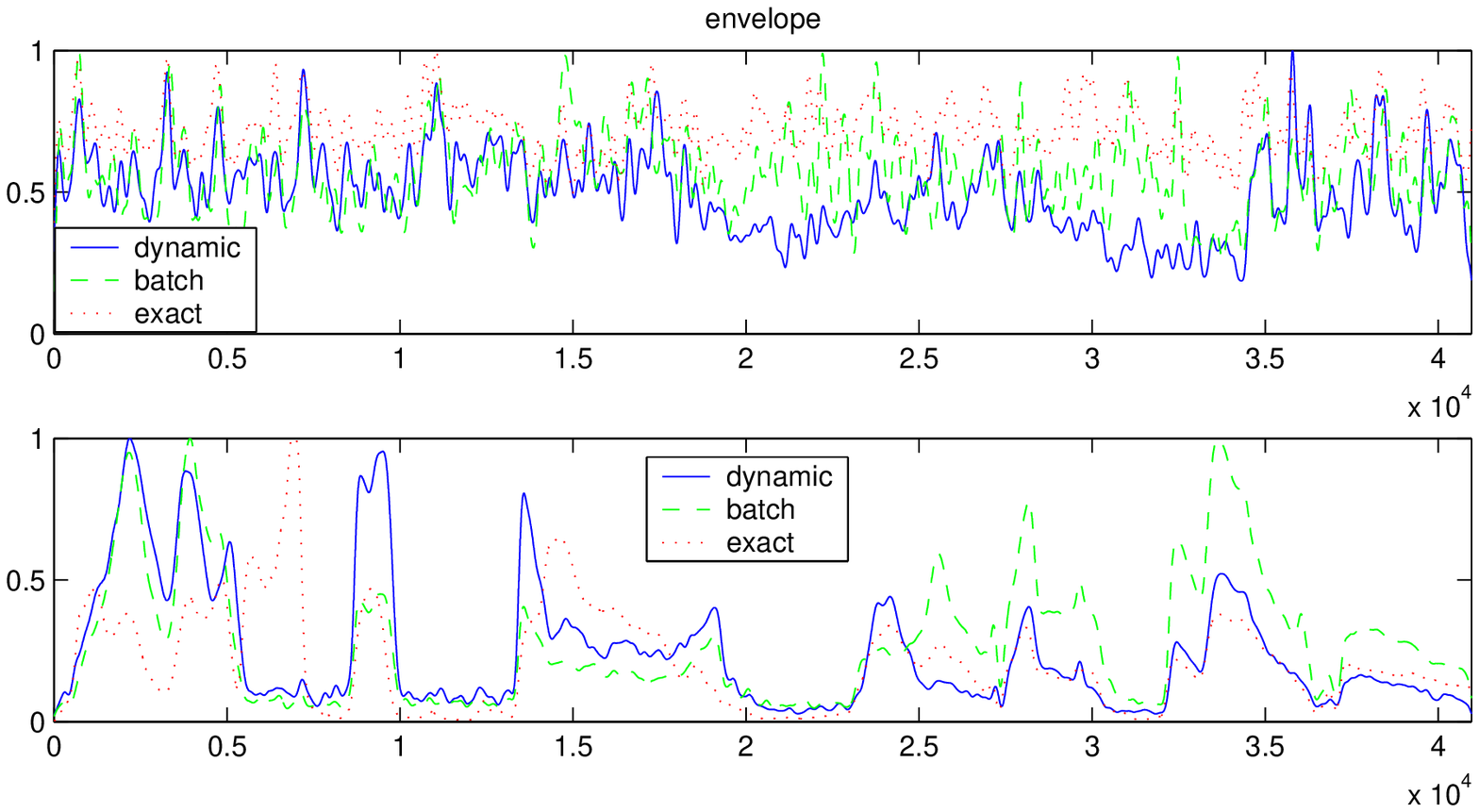}}
\caption{} \label{f-n.2}
\end{figure}

\end{document}